\documentclass[runningheads]{llncs}

\usepackage[utf8]{inputenc} 
\usepackage[T1]{fontenc}    
\usepackage{hyperref}       
\usepackage{url}            
\usepackage{booktabs}       
\usepackage{amsfonts}       
\usepackage{nicefrac}       
\usepackage{microtype}      

\usepackage{doi}

\usepackage{amsmath, amssymb, bm}

\usepackage{xcolor} 

\usepackage{threeparttable} 

\usepackage{graphicx}
\makeatletter
\newcommand{\printfnsymbol}[1]{%
  \textsuperscript{\@fnsymbol{#1}}%
}
\makeatother

\usepackage[title]{appendix}

\begin{document}
\title{Dominant Eigenvalue-Eigenvector Pair\\
Estimation via Graph Infection}
%
%
\author{Kaiyuan Yang\inst{1} \thanks{Equal Contribution.} \and
Li Xia\inst{2} \printfnsymbol{1} \and
Y.C. Tay\inst{3}}
\authorrunning{K. Yang and L. Xia et al.}
%
\institute{Department of Quantitative Biomedicine,\\ University of Zurich,
8057 Zurich, Switzerland
\and
Department of Statistics and Data Science,\\
National University of Singapore, Singapore 117546
\and
Department of Computer Science,\\
National University of Singapore, Singapore 117417\\
\email{dcstayyc@nus.edu.sg}}
\maketitle              
\begin{abstract}
    We present a novel method to estimate the dominant eigenvalue and eigenvector pair
    of any non-negative real matrix via graph infection.
    The key idea in our technique lies in approximating the solution to the first-order matrix ordinary differential equation (ODE) with the Euler method.
    Graphs, which can be weighted, directed, and with loops, are first converted to its adjacency matrix A.
    Then by a naive infection model for graphs, we establish the corresponding first-order matrix ODE,
    through which A's dominant eigenvalue is revealed by the fastest growing term.
    When there are multiple dominant eigenvalues of the same magnitude,
    the classical power iteration method can fail.
    In contrast, our method can converge to the dominant eigenvalue even when same-magnitude counterparts exist, be it complex or opposite in sign. 
    We conduct several experiments comparing the convergence between our method and power iteration.
    Our results show clear advantages over power iteration for tree graphs, bipartite graphs, directed graphs with periods, and Markov chains with spider-traps.
    To our knowledge, this is the first work that estimates dominant eigenvalue and eigenvector pair
    from the perspective of a dynamical system and matrix ODE.
    We believe our method can be adopted as an alternative to power iteration, especially for graphs.

\keywords{Dominant Eigenvalue \and Eigenvector Centrality \and First-order Matrix ODE \and Euler Method \and Graph Infection}
\end{abstract}

\section{Introduction}





Graph epidemic models seek to describe the dynamics of contagious disease transmission over networks~\cite{kiss_EoN_textbook_Epidemic_2017mathematics,st_andrew_dobson_epidemic_book}.
The infectious disease transmits from a node to its neighbors via connecting edges over the network.
Spread of the epidemic is affected by multiple factors such as
the infection rate,
the recovery duration,
and infection severity,
and particularly to graphs, the network topology and the mobility of the network structure.
Graph epidemic models can encode richer and more sophisticated architectures than traditional compartmental epidemic models~\cite{kiss_EoN_textbook_Epidemic_2017mathematics,st_andrew_dobson_epidemic_book}.

For disease spread on networks,
the principal eigenvalue of the network's adjacency matrix has long been shown to be an important factor on the dynamics of disease~\cite{wang2003epidemic_CMU,ganesh2005effect_eigenNetwork}.
In fact, the well-known basic reproduction number `R0' is itself the dominant eigenvalue of the next generation matrix~\cite{diekmann1990_R0_definition,diekmann2010_R0_construction}.
Interestingly, despite the clear connection between the dominant eigenvalue and the network epidemics,
few have attempted to approach the principal eigenvalue in the reverse manner.
Here we seek to answer the question:
Can we elucidate the principal eigenvalue from the progression of the infection spread over the associated network?

Another motivation of ours comes from the extant issues on eigenvalue computation from the classical power iteration method.
Power iteration method~\cite{golub_origin_powermethod_2000}, or power method, has been widely used for computation of the principal eigenvalue.
However, when multiple dominant eigenvlaues of equal modulus exist, power iteration method cannot converge~\cite{quarteroni2010numerical_powermethodpage192}.
Note that such failed convergence of power method is not uncommon, especially for graphs.
One prominent failure is for any bipartite graph:
all non-zero eigenvalues of a bipartite graph's adjacency matrix come in real number pairs that are of opposite sign.
Thus if there exists a dominant eigenvalue $\lambda$ for the adjacency matrix of the bipartite graph, then so is $-\lambda$.
Power iteration method can also fail in Markov chains with periods.
For example, there are as many eigenvalues equally spaced around the unit circle as the period of the periodic unichain~\cite{gallager2013stochastic_unichain_circle}.
Unless the Markov chain is an ergodic unichain which has only one dominant eigenvalue of 1,
problematic convergence to the dominant eigenvector from power iteration should be expected.

In this paper, we aim to develop an epidemic-based method to estimate the principal eigenvalue and eigenvector of a network, and compare its applications with the classical power iteration method.
Our perspective of using the fastest growing term in general solution of graph infection ordinary differential equation (ODE) to `reverse engineer' the principal eigenvalue is new.
Our proposed alternative method can also overcome several limitations of power iteration method.
In particular, our method works better when there are multiple dominant eigenvalues of opposite sign or are complex conjugates.















\section{Our Method: Inspired by Euler Method for Graph Infection ODE}

In this section we present our proposed method.
See Table~\ref{tab:table_notations} for the notations and symbols.
\begin{table}
	\caption{Table of Notations}
	\centering
	\begin{tabular}{ll}
		\toprule
		Notation     & Description \\
		\midrule
		$\mathcal{G}$ & Graph, $\mathcal{G}=\{\mathcal{V}, \mathcal{E}\}$. Graphs can be \textbf{weighted, directed, and with loops} \\
		$\mathcal{V}$ & Set of vertices of the graph, $\mathcal{V} = \{v_1,..., v_N\}$, $\lvert \mathcal{V} \rvert = N$\\
		$\mathcal{E}$ & Set of edges of the graph \\
		$\mathcal{N}(i)$ & Set of in-degree neighbour nodes of node $v_i$ \\
        $\bm{A}$ & Adjacency matrix of the graph, $\bm{A}\in \mathbb{R}_{\geq0}^{N \times N}$\\
        $\beta$ & Ratio of infection severity transmitted per unit time, $\beta \in \mathbb{R}_{>0}$ \\
        $x_i(t)$ & Infection severity of node $v_i$ at time $t$, $x_i(t) \in \mathbb{R}_{\geq0}$\\
        $\mathbf{x}(t)$ & Vector form of the infection severity of each node at time $t$\\
        $I(t)$ & Total infection of the graph at time $t$, $I(t)=\sum_{i=1}^N x_i(t)$\\
        $\lambda_k$ & $k$-th eigenvalue of $\bm{A}$, $\lambda_k \in \mathbb{C}$. $|\lambda_1| \geq |\lambda_2| \geq ... \geq |\lambda_N|$\\
        $\bm{\mu}_k$ & $k$-th eigenvector of $\bm{A}$, $\bm{A}\bm{\mu}_k=\lambda_k \bm{\mu}_k$ and $\bm{\mu}_k \neq \mathbf{0}$\\
        $r_k$ & Algebraic multiplicity of the $k$-th eigenvalue of $\bm{A}$\\
        $p_k$ & Geometric multiplicity of the $k$-th eigenvalue of $\bm{A}$\\
        $d_k$ & Deficit of the $k$-th eigenvalue of $\bm{A}$, $d_k = r_k - p_k$\\
        $C_k$ & Coefficient of $\bm{\mu}_k$ in the linear combination of eigenvectors for initial condition $\mathbf{x}(0)$\\
        $m$ & Slope of the \textbf{secant line} with step-size $\Delta t$ for log-scale plot $I(t)$ vs time\\
		\bottomrule
	\end{tabular}
	\label{tab:table_notations}
\end{table}
The key idea in our technique lies in solving the matrix ODE that describes the graph infection process.
The fastest growing term of the general solution to the ODE will reveal the adjacency matrix's dominant eigenvalue.
This idea comes from a naive infection model for a graph,
\textbf{but note that our estimation technique is numerical and does not simulate infection}.

\subsection{Graph Infection and Adjacency Matrix}

For our setup, \textbf{our graph is static}, i.e., the network structure does not vary as the infection process unfolds.
For a given graph $\mathcal{G}$, the edges $\mathcal{E}$ connecting nodes $\mathcal{V}$ over the network do not rewire or change weights once the infection starts.
The graph can be \textbf{weighted, directed, and with loops}.
For such a graph network, the adjacency matrix $\bm{A}$ describes the connections between its nodes,
in which each entry $a_{ij}$ represents the edge from node $j$ to node $i$.
Matrix $\bm{A}$ can be symmetric for undirected networks, and potentially asymmetric for directed networks.
Elements of $\bm{A}$ can be zero and one for unweighted connections.
Or in the case of weighted connections, $\bm{A}$ may be an arbitrary non-negative matrix.
Thus our $\bm{A}\in \mathbb{R}_{\geq0}^{N \times N}$, where $N$ is the number of nodes.

\begin{figure}[!ht]
	\centering
	\includegraphics[width=0.55\textwidth]{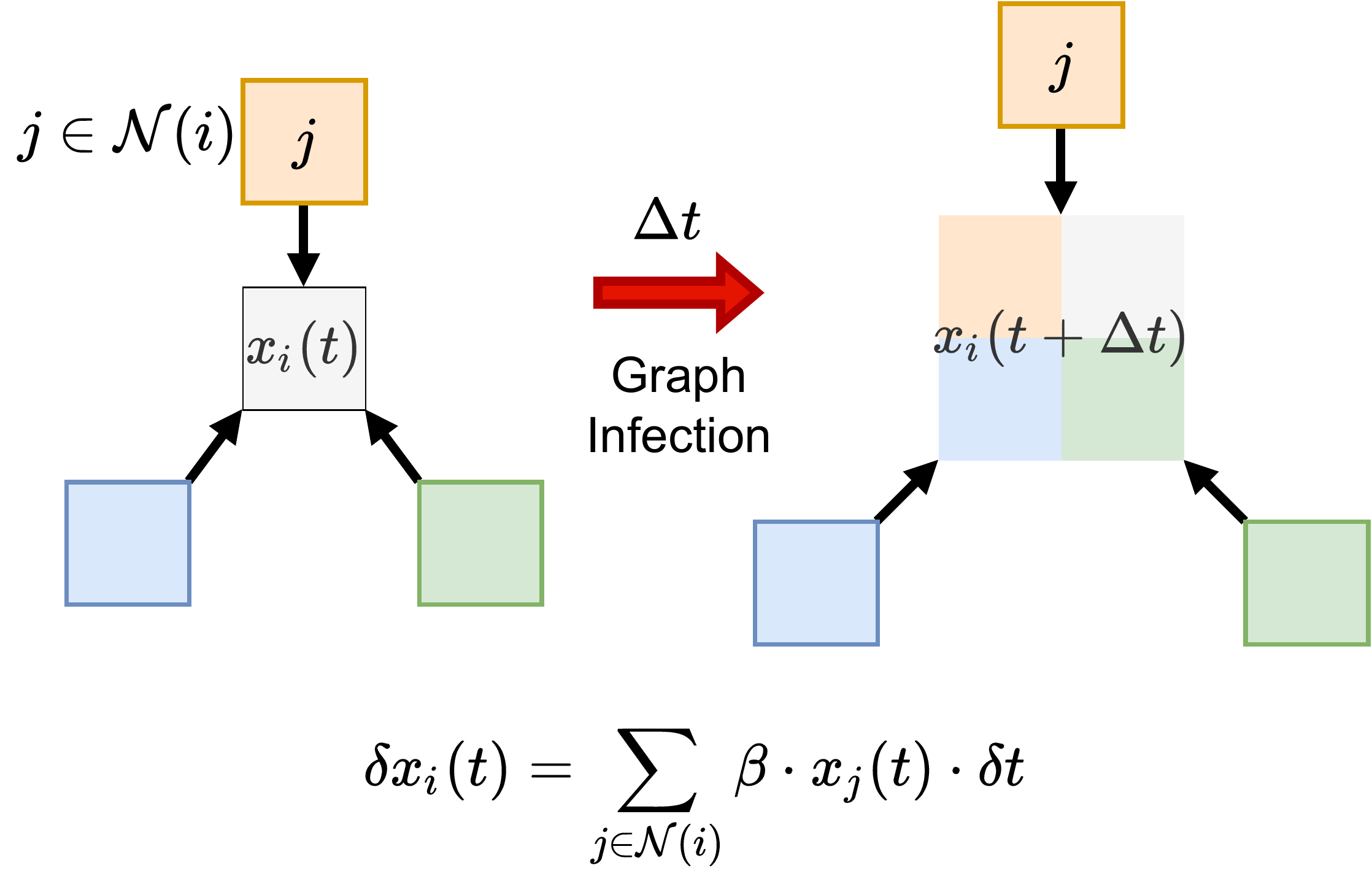}
	\caption{Graph infection with infection rate of $\beta$ and the difference equation of the infection severity $x_i(t)$ at node $v_i$.
	Infection spreads from infected neighbouring nodes $\mathcal{N}(i)$ via in-degree edges. Infection severity accumulates, think of virus count in a host.}
	\label{fig:fig_I_vs_dt}
\end{figure}

A basic example demonstrating the mechanism of graph infection for one node is shown in Figure~\ref{fig:fig_I_vs_dt}.
The rate of change of the infection severity $x_i(t)$ at node $v_i$ can be intuitively expressed as a difference equation shown in the figure.
In order to describe the infection change for the graph, we turn to matrix ODE.

\subsection{General Solution to Graph Infection ODE}

\begin{figure}
	\centering
	\includegraphics[width=\textwidth]{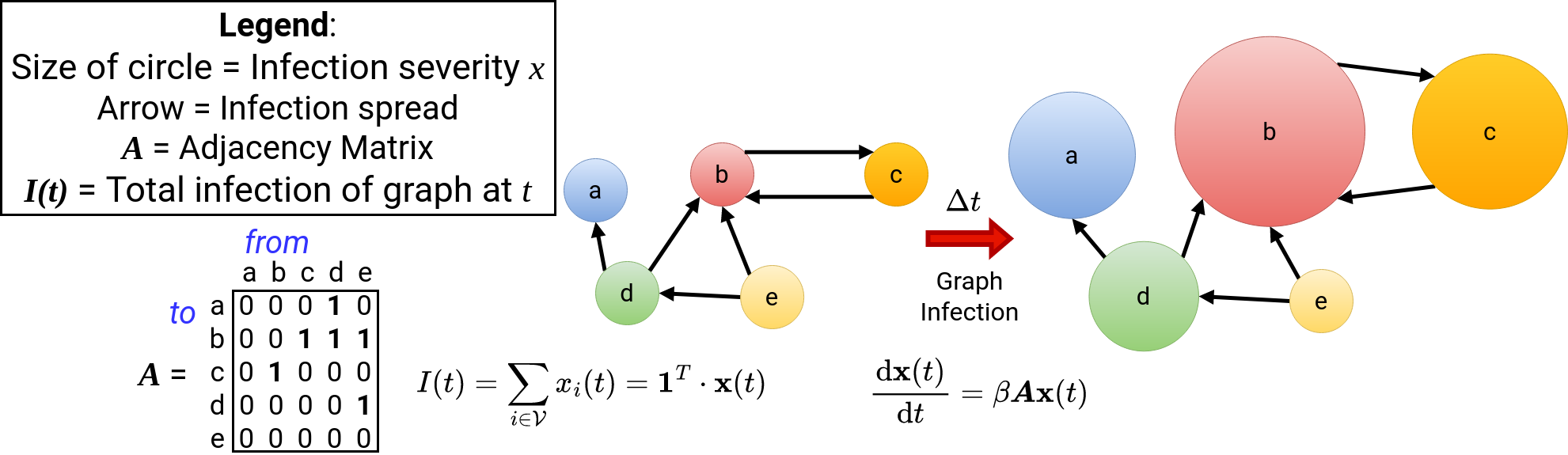}
	\caption{An example network and its evolution of infection severity.}
	\label{fig:fig_Ax}
\end{figure}

\textbf{The setup for our numerical method is based on the change of severity of infection for each node in the network over time.}
Conceptually, we can view this model as a description of infection
severity of nodes, where nodes
with severely infected neighbors would receive more severe transmission.
See Figure~\ref{fig:fig_Ax} for a sketch of the infection severity ODE idea.
Alternatively, one may perceive this model as a pest which
continuously multiplies and invades neighboring nodes over a network.
\textbf{Note that we make the evolution of severity deterministic,
and thus remove stochastic elements in the model.}

The severity of node $i$ at time $t$, $x_i(t)$,
is related to its neighboring nodes' severity.
$x_i(t)$ changes with time according to this differential equation:
\begin{equation}
  \frac{\mathrm{d} x_i(t) }{\mathrm{d} t} = \beta \sum_{j=1}^{N} \bm{A}_{ij} x_j(t)
  \label{eq:Euler_dzdt}
\end{equation}
Using $\mathbf{x}(t)$ to represent the column vector of each node's infection severity at time $t$,
Equation~\ref{eq:Euler_dzdt} becomes:
\begin{align}
  \frac{\mathrm{d} \mathbf{x}(t)}{\mathrm{d} t} &= \beta \bm{A}\mathbf{x}(t)
  \label{eq:Euler_dzdt_Az}
\end{align}

\textbf{Note the general solution to the above matrix ODE of Equation~\ref{eq:Euler_dzdt_Az} can be divided into \textit{four} cases based on the eigenvalues of $\bm{A}$.}
Since the characteristic equation of $\bm{A}$ is an $N$th-order polynomial with real coefficients,
$\bm{A}$ has exactly $N$ eigenvalues (including repetitions if any)
that are either a real number or a pair of complex conjugate numbers.
Here we distinguish the four cases based on whether the eigenvalue is real or complex and on whether it is repeated.

\subsubsection{Case One: Real Distinct Eigenvalues}

The first case is when the adjacency matrix $\bm{A}$ has $N$ real and distinct eigenvalues.
Then we can decompose the initial condition into $N$ linearly independent eigenvectors,
\begin{equation}
    \mathbf{x}(0) = \sum_{k=1}^{N} C_k \bm{\mu}_k \nonumber
    \label{eq:linear_decompose_x0}
\end{equation}
In this case, the matrix ODE Equation~\ref{eq:Euler_dzdt_Az} has general solution:
\begin{equation}
    \mathbf{x}(t) = \sum_{k=1}^{N} C_k e^{\beta \lambda_k t}\bm{\mu}_k
    \label{eq:Euler_dzdt_Az_genSol}
\end{equation}
where $\lambda_1, \lambda_2, ..., \lambda_N$ are the eigenvalues of the adjacency matrix $\bm{A}$,
and $\bm{\mu}_1, \bm{\mu}_2, ..., \bm{\mu}_N$ are the corresponding eigenvectors.


\subsubsection{Case Two: Complex Distinct Eigenvalues}

In the second case, $\bm{A}$ still has $N$ distinct eigenvalues,
however, some of the eigenvalues are complex.
Since $\bm{A}$ is real, if a complex number $\lambda_k = a+ib$ is an eigenvalue of our $\bm{A}$,
then its complex conjugate $\overline{\lambda_k} = a-ib$ is also an eigenvalue, with $a, b \in \mathbb{R}$.

Let $\lambda_k = a+ib$ be a complex eigenvalue of $\bm{A}$, and $\bm{\mu}_k = \mathbf{a} + i\mathbf{b}$ be the corresponding eigenvector where vectors $\mathbf{a}$ and $\mathbf{b}$ have only real entries.
Since $\bm{A}$ still has $N$ linearly independent eigenvectors,
the general solution to the matrix ODE Equation~\ref{eq:Euler_dzdt_Az} remains in the form of Equation~\ref{eq:Euler_dzdt_Az_genSol}.
However, the general solution will contain the following two terms in place of the conjugate complex eigenvalues $\lambda_k$ and $\overline{\lambda_k}$ and their eigenvectors:
\begin{equation}
    C_p e^{\beta a t} \left( \mathbf{a} \cos(\beta bt) - \mathbf{b} \sin(\beta bt) \right) + C_q e^{\beta a t} \left( \mathbf{a} \sin(\beta bt) + \mathbf{b} \cos(\beta bt) \right)
    \label{eq:complex_gen_sol}
\end{equation}
where $C_p$ and $C_q$ are some scalars determined by the initial condition $\mathbf{x}(0)$.

\subsubsection{Case Three: Repeated Complete Eigenvalues}

When the characteristic equation of $\bm{A}$ has a multiple root,
the associated eigenvalue is called a repeated eigenvalue.
Suppose $\lambda_k$ is a real multiple-root of the characteristic equation of $\bm{A}$,
then the algebraic multiplicity of $\lambda_k$,
defined as $r_k$,
is the number of times $\lambda_k$ appears as a root.
On the other hand, geometric multiplicity of $\lambda_k$,
defined as $p_k$,
is the number of linearly independent eigenvectors associated with $\lambda_k$.
Here we distinguish whether the algebraic multiplicity and geometric multiplicity are equal or not.

In Case Three, the algebraic multiplicity and geometric multiplicity for the repeated eigenvalue are equal,
or $r_k = p_k$.
And we say the repeated eigenvalue is \textbf{complete}, as in without deficit.
Let $\lambda_k$ be the repeated eigenvalue with geometric multiplicity $p_k$, where $r_k = p_k$,
and $\bm{\mu}_k, \mathbf{w}_{1}, \mathbf{w}_{2}, ..., \mathbf{w}_{p_k-1}$ are the $p_k$ linearly independent ordinary eigenvectors associated with $\lambda_k$.
Thus despite repeated eigenvalues in $\bm{A}$,
there are still $N$ linearly independent eigenvectors.
The general solution to the matrix ODE Equation~\ref{eq:Euler_dzdt_Az} remains in the form of Equation~\ref{eq:Euler_dzdt_Az_genSol},
but contains the following term for the repeated complete eigenvalue $\lambda_k$:
\begin{equation}
    C_1 e^{\beta \lambda_k t} \bm{\mu}_k + C_2 e^{\beta \lambda_k t} \mathbf{w}_1 + ... + C_{p_k} e^{\beta \lambda_k t} \mathbf{w}_{p_k-1}
    \label{eq:fastest_term_complete_repeat}
\end{equation}
where $C_1$ to $C_{p_k}$ are some scalars determined by the initial condition $\mathbf{x}(0)$.

\subsubsection{Case Four: Repeated Defective Eigenvalues}

In Case Four, there exists a repeated eigenvalue $\lambda_k$ for which the geometric multiplicity is smaller than the algebraic multiplicity,
or $r_k > p_k$.
Hence, there are fewer than $N$ linearly independent eigenvectors for $\bm{A}$.
We call such a repeated eigenvalue $\lambda_k$ as \textbf{defective},
and $d_k = r_k - p_k$ is the \textbf{deficit} of $\lambda_k$.

In order to compensate for the defect gap in the multiplicity,
we need to construct $d_k$ number of \textit{generalized eigenvectors} for the ODE general solution.
These generalized eigenvectors are associated with the defective repeated eigenvalue $\lambda_k$
and based on the $p_k$ linearly independent ordinary eigenvectors.

Let $\lambda_k$ be the repeated eigenvalue with deficit of $d_k$,
and $\bm{\mu}_k, \mathbf{w}_{1}, \mathbf{w}_{2}, ..., \mathbf{w}_{p_k-1}$ are the $p_k$ linearly independent ordinary eigenvectors associated with $\lambda_k$.
Note that the eigenspace must have dimension of at least one.
Therefore $\bm{\mu}_k$ is guaranteed to be available for us to construct the generalized eigenvectors.
In the simplest two scenarios, we can have either one generalized eigenvector to find or we have only one ordinary eigenvector available.

In the first scenario of only needing to find one generalized eigenvector,
or when we have $d_k = 1$ and $r_k  = p_k + 1$,
the general solution to the matrix ODE Equation~\ref{eq:Euler_dzdt_Az}
will have the following term for the defective repeated eigenvalue $\lambda_k$:
\begin{equation}
    C_1 e^{\beta \lambda_k t} \bm{\mu}_k + C_2 e^{\beta \lambda_k t} \mathbf{w}_1 + ... + C_{p_k} e^{\beta \lambda_k t} \mathbf{w}_{p_k-1} + C_{p_k + 1} e^{\beta \lambda_k t} \left( t \bm{\mu}_k + \mathbf{g} \right)
    \label{eq:defect_deficit_is_1}
\end{equation}
where $C_1$ to $C_{p_k+1}$ are some scalars determined by the initial condition $\mathbf{x}(0)$,
and vector $\mathbf{g}$ is a generalized eigenvector
such that $(\bm{A} - \lambda_k \bm{I})\mathbf{g} = \bm{\mu}_k$.

In the second scenario of having only one ordinary eigenvector $\bm{\mu}_k$ associated with $\lambda_k$,
or when $p_k = 1$ and $r_k = d_k + 1$,
the general solution to the matrix ODE Equation~\ref{eq:Euler_dzdt_Az}
will have the following term for the defective repeated eigenvalue $\lambda_k$:
\begin{equation}\label{eq:defect_geo_is_1}
\begin{split}
    C_1 e^{\beta \lambda_k t} \bm{\mu}_k &+ C_2 e^{\beta \lambda_k t} \left( t \bm{\mu}_k + \mathbf{\tilde{g}}_1 \right) + C_3 e^{\beta \lambda_k t} \left( \frac{t^2}{2!}  \bm{\mu}_k + t \mathbf{\tilde{g}}_1 + \mathbf{\tilde{g}}_2 \right) + ... \\
    &{...}+ C_{d_k + 1} e^{\beta \lambda_k t} \left( \frac{t^{d_k}}{d_k!}  \bm{\mu}_k + \frac{t^{d_k-1}}{(d_k-1)!} \mathbf{\tilde{g}}_1 + ... + \frac{t^2}{2!} \mathbf{\tilde{g}}_{d_k-2} + t \mathbf{\tilde{g}}_{d_k-1} +  \mathbf{\tilde{g}}_{d_k} \right)
\end{split}
\end{equation}
where $C_1$ to $C_{d_k+1}$ are some scalars determined by the initial condition $\mathbf{x}(0)$,
and vectors $\mathbf{\tilde{g}}_{1}$ to $\mathbf{\tilde{g}}_{d_k}$ are generalized eigenvectors
such that
$(\bm{A} - \lambda_k \bm{I})\mathbf{\tilde{g}}_1 = \bm{\mu}_k$
and
$(\bm{A} - \lambda_k \bm{I})\mathbf{\tilde{g}}_{j+1} = \mathbf{\tilde{g}}_{j}$ for $j\in \{1, ..., d_k\}$.

In between the two extreme scenarios,
when the deficit $1<d_k<r_k-1$,
the chains of generalized eigenvectors can be made up of some combinations of the $p_k$ ordinary eigenvectors associated with the defective repeated eigenvalue $\lambda_k$.
Although the structure of the chains of generalized eignvectors can be complicated,
we will still end up with $r_k$ number of linearly independent solution vectors involving chains of generalized eigenvectors.
Therefore the general solution to the matrix ODE Equation~\ref{eq:Euler_dzdt_Az} can contain $r_k$ number of terms of the ordinary and generalized eigenvectors arranged in chains,
resembling the terms of Equation~\ref{eq:defect_deficit_is_1} and Equation~\ref{eq:defect_geo_is_1}.

Note that it is possible for $\bm{A}$ to have repeated complex conjugate pair of eigenvalues.
The method involving generalized eigenvectors as discussed above works for defective complex eigenvalue as well.
For more details on the multiple eigenvalue solutions to matrix ODE,
please refer to Chapter 5 in~\cite{edwards2008differential}.
However,
the repeated complex eigenvalues are not a major concern for us because of the Perron-Frobenius theorem
as we will discuss in the following section.

\subsection{Fastest Growing Term and Perron-Frobenius Theorem}

We first examine the graph infection ODE's general solution from the perspective of the long-term behavior of the network.
According to the \textbf{Perron–Frobenius theorem} extended to \textbf{non-negative real matrices},
there exists a real non-negative eigenvalue
that is greater than or equal to
the absolute values of all other eigenvalues of that matrix.
So
\begin{align*}
    \exists \lambda_1 .\; \lambda_1 \in \mathbb{R}_{\geq0}\\
    \forall \lambda_k .\; k > 1 \Rightarrow \lambda_1 \geq  |\lambda_k|
\end{align*}

In the Case One of real and distinct eigenvalues,
the dominance of the principal real eigenvalue will mean the following:
\begin{align*}
    \lambda_1 > \lambda_k \Rightarrow e^{\lambda_1 t} > e^{\lambda_k t}
\end{align*}
Thus $\beta\lambda_1 t$ becomes the dominant exponent term for Equation~\ref{eq:Euler_dzdt_Az_genSol}.
\textbf{Since it is exponential, $\exp \left( \beta \lambda_1 t \right)$ very quickly overwhelms all the other terms in the summation
in the Case One},
as the eigenvalues are real and distinct.
The dynamics along the principal eigenvalue will dominate the long-term behavior in Equation~\ref{eq:Euler_dzdt_Az}.
Therefore the fastest growing term in the general solution for Equation~\ref{eq:Euler_dzdt_Az} corresponds to:
\begin{align}
    \mathbf{x}(t) &\approx C_1 e^{\beta \lambda_1 t} \bm{\mu}_1
    \label{eq:fastest_growing_zt}
\end{align}
i.e., $\lambda_1$ and $\bm{\mu}_1$ are revealed through the change in $\mathbf{x}(t)$.

\textbf{Note that this dominant exponent also works for complex eigenvalues in the Case Two of having complex eigenvalues.}
Since for $k > 1$, $\lambda_k \in \mathbb{C}$,
let $\lambda_k = a + ib$, for some real and non-zero numbers $a$ and $b$.
Then
\begin{equation*}
    \lambda_1 \geq  |\lambda_k| = |a + ib| \Rightarrow
    \lambda_1 > a \Rightarrow
    e^{\lambda_1 t} > e^{a t}
\end{equation*}
Hence, the exponent with principal real eigenvalue will still dominate the long-term behavior of the general solution to the matrix ODE,
even when there are complex eigenvalues with the same norm as the principal eigenvalue in Equation~\ref{eq:complex_gen_sol}.
Thus we can derive $\lambda_1$ and $\bm{\mu}_1$ using Equation~\ref{eq:fastest_growing_zt} similar to Case One.

When there are repeated eigenvalues as in Case Three and Case Four,
unless the principal real eigenvalue has algebraic multiplicity greater than $1$,
the asymptotic behavior of the ODE's general solution is still Equation~\ref{eq:fastest_growing_zt}.
Now suppose the repeated eigenvalue is the principal real eigenvalue $\lambda_1$ and $r_1 > 1$.
In Case Three,
since the algebraic multiplicity and geometric multiplicity are the same,
$r_1 = p_1$,
the fastest growing term in the general solution to ODE is Equation~\ref{eq:fastest_term_complete_repeat},
which is a scalar multiple of $\exp \left( \beta \lambda_1 t \right)$,
and converges to a vector spanned by the $p_k$ linearly independent ordinary eigenvectors associated with $\lambda_1$.
Similarly for Case Four,
the fastest growing term in the general solution to ODE is a scalar multiple of $\exp \left( \beta \lambda_1 t \right)$.
However, the associated vector in Case Four is not convergent.

Table~\ref{tab:table_fastestGrowingTerm} gives a summary of the fastest growing term to the ODE for Case One to Four.
\begin{table}
\begin{threeparttable}[b]
	\caption{Fastest Growing Term in the general solution to ODE Equation~\ref{eq:Euler_dzdt_Az}}
	\centering
	\begin{tabular}{l|c|c}
		\toprule
		Case     & Dominant Eigenvalues $|\lambda_k|  = \lambda_1$    & Fastest Growing Term \\
		\midrule
		1: Real Distinct Eigenvalues & $-\lambda_k = \lambda_1.\; \lambda_k \in \mathbb{R}_{\leq0}$  & $C_1 e^{\beta \lambda_1 t} \bm{\mu}_1$  \\
		2: Complex Distinct Eigenvalues & $|\lambda_k|  = \lambda_1.\; \lambda_k \in \mathbb{C}$ & $C_1 e^{\beta \lambda_1 t} \bm{\mu}_1$ \\
		3: Repeated Complete Eigenvalues & $\lambda_k = \lambda_1.\; \lambda_k \in \mathbb{R}_{\geq0}$  & $e^{\beta \lambda_1 t} \left(C_1 \bm{\mu}_1 + C_k\bm{\mu}_k \right)$  \\
		4: Repeated Defective Eigenvalues & $\lambda_k = \lambda_1.\; \lambda_k \in \mathbb{R}_{\geq0}$  & $e^{\beta \lambda_1 t} \left( C_1 \bm{\mu}_1 + C_k \left( t \bm{\mu}_1 + \mathbf{g} \right) \right)$\tnote{*}  \\
		\bottomrule
	\end{tabular}
     \begin{tablenotes}
       \item [*] For Case Four, here we only show the simplest scenario of $d_1 = 1, r_1 = 2$, and $\mathbf{g}$ is a generalized eigenvector such that $\bm{A} \mathbf{g} = \lambda_1 \mathbf{g} + \bm{\mu}_1$.
     \end{tablenotes}
	\label{tab:table_fastestGrowingTerm}
\end{threeparttable}
\end{table}
Suppose there are multiple eigenvalues of the same magnitude or modulus,
say $\lambda_k$ and $\lambda_1$.
As shown in Table~\ref{tab:table_fastestGrowingTerm},
in all four Cases,
the fastest growing term converges to
the Perron-Frobenius real dominant eigenvalue $\lambda_1 \in \mathbb{R}_{\geq0}$.
In Case One and Two,
despite the presence of a $\lambda_k$ of opposite sign or is a complex conjugate,
the fastest growing term is convergent to both the $\lambda_1$ and $\bm{\mu}_1$.
As we shall see later,
this marks the major advantage over the classical power iteration method.

\subsection{Euler Method to Approximate Matrix ODE}

Let $I(t)$ be the total severity of infection for the graph at time $t$:
\begin{align*}
    I(t) &= \sum_{i=1}^{N} x_i(t) = \bm{1}^T \cdot \mathbf{x}(t)
\end{align*}

We turn to \textbf{Euler method}.
From matrix differential Equation~\ref{eq:Euler_dzdt_Az}, we can estimate change of $\mathbf{x}(t)$ by the \textbf{finite difference approximation}:
\begin{align}
  \mathbf{x}(t + \Delta t) &\approx \mathbf{x}(t) + \frac{\mathrm{d} \mathbf{x}(t)}{\mathrm{d} t}\Delta t \nonumber \\
  &= \mathbf{x}(t) + \beta \bm{A}\mathbf{x}(t)\Delta t \quad \text{by Equation~\ref{eq:Euler_dzdt_Az}}
  \label{eq:reveal_eigVec}
\end{align}

Then we rewrite the finite difference for severity of graph as:
\begin{align}
  I(t+\Delta t) &= \bm{1}^T \cdot \mathbf{x}(t + \Delta t) \nonumber \\
  &\approx \bm{1}^T \cdot \left( \mathbf{x}(t) + \beta \bm{A}\mathbf{x}(t)\Delta t \right) \nonumber \\
  &\approx \bm{1}^T \cdot C_1 e^{\beta \lambda_1 t} \bm{\mu}_1 \cdot \left( 1 + \beta \lambda_1 \Delta t \right) \nonumber \nonumber \\
  &= \bm{1}^T \cdot \mathbf{x}(t) \cdot \left( 1 + \beta \lambda_1 \Delta t \right) \nonumber \\
  &= I(t) \cdot \left( 1 + \beta \lambda_1 \Delta t \right) \label{eq:ItdtvsIt}
\end{align}
Note that here we substitute $\mathbf{x}(t)$ with the fastest growing term for Case One and Case Two
from Table~\ref{tab:table_fastestGrowingTerm}.
But it should be obvious that when there are repeated and complete dominant eigenvalues,
the fastest growing term of Case Three from Table~\ref{tab:table_fastestGrowingTerm} substituting $\mathbf{x}(t)$ will also give the same derivation here at Equation~\ref{eq:ItdtvsIt}.

\subsection{Secant Line of Euler Method}

\begin{figure}[!b]
	\centering
	\includegraphics[width=\textwidth]{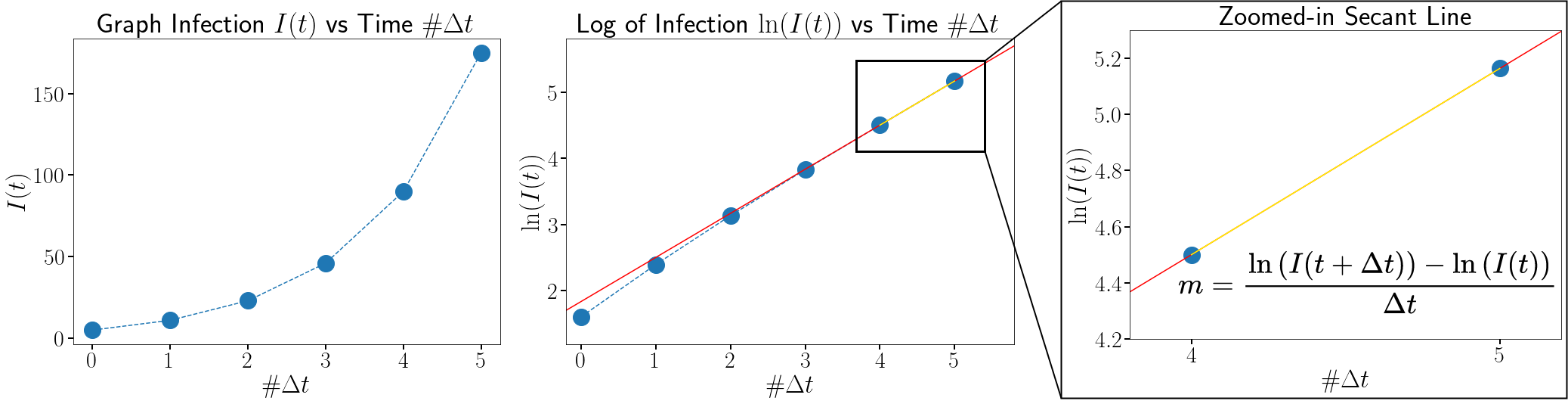}
	\caption{Secant line of graph infection plot using the adjacency matrix $\bm{A}$ from the network in Figure~\ref{fig:fig_Ax}, with $\mathbf{x}(0) = \bm{1}$, $\beta = 1$, $\Delta t = 1$.
	\textbf{(Left)} Plot of graph infection $I(t)$ vs time in number of discrete time-steps $\#\Delta t$.
	Blue dotted lines connecting the data points are the secant lines.
	\textbf{(Middle)} Plot the log of graph infection $\ln (I(t))$ vs discrete time-steps.
	Yellow segment is a highlighted secant line connecting two graph infection measurements.
	Red line is the slope $m$ of the yellow secant line.
	\textbf{(Right)} Zoomed-in view of the secant line. The slope $m$ is used to estimate the dominant eigenvalue and eigenvector.}
	\label{fig:fig_secant_explain}
\end{figure}

Next we make use of the finite difference of the graph infection severity with secant line. Figure~\ref{fig:fig_secant_explain} shows a graphical explanation of the role of secant line and Euler method in our approximation.
We plot the \textbf{log of total severity} of the graph at time $t$, $\ln(I(t))$, versus time in number of discrete time-steps,
as shown in Figure~\ref{fig:fig_secant_explain}.
Take two points along the $\ln(I(t))$ curve separated by $\Delta t$,
and \textbf{draw a secant line},
that secant line will have slope $m$ defined by:
\begin{align*}
  m &= \frac{\ln \left( I(t+\Delta t) \right) - \ln \left( I(t) \right)}{\Delta t} \\
  &= \frac{\ln \left(
   1 + \beta \lambda_1 \Delta t
  \right)}{\Delta t} \quad \text{by Equation~\ref{eq:ItdtvsIt}}
\end{align*}
Note that the slope $m$ of the secant line can be measured for a given time-step $\Delta t$ as shown in Figure~\ref{fig:fig_secant_explain}.
Therefore the slope $m$ of the secant line and the true dominant eigenvalue $\lambda_1$ is related by:
\begin{align}
  e^{m\Delta t} &= 1 + \beta \lambda_1\Delta t \nonumber \\
  \therefore \lambda_1 &= \frac{e^{m\Delta t} - 1}{\beta \Delta t}
  \label{eq:last_dance}
\end{align}
The associated dominant eigenvector $\bm{\mu}_1$ is revealed by Equation~\ref{eq:reveal_eigVec} along with the iterative process.
Since our method depends on the long-term behavior of the dynamical system,
we iteratively calculate the secant line slope $m$ after each time step is taken.
After sufficient number of time steps,
we start to zoom into large $t$ when the matrix ODE solution is dominated by the fastest growing term.
This secant-line-based method thus allows us to obtain a fair estimation of the dominant eigenvalue $\lambda_1$ via measuring $m$ for a chosen $\beta$ and $\Delta t$ as elucidated in Equation~\ref{eq:last_dance}.
\textit{(Note that infection rate $\beta$ is introduced to help set up the context and perspective of epidemic and infection adopted by our method.
However, $\beta$ can be assumed to be always $1$ for simplicity,
i.e. $100\%$ of a node's infection severity will be transmitted through a connecting outward edge per unit time.)}


\section{Theoretical Comparison with Power Iteration}

Table~\ref{tab:table_compare_PowerMethod} shows the comparison between our method with the classical power iteration method.
\begin{table}
\begin{threeparttable}[b]
	\caption{Convergence comparison: Our method with the power iteration method}
	\centering
	\begin{tabular}{llcccc}
		\toprule
		 & & \multicolumn{2}{c}{Converge to $\lambda_1$?} & \multicolumn{2}{c}{Converge to $\bm{\mu}_1$?} \\
		\cmidrule(r){3-4} \cmidrule(r){5-6}
		Case & Dominant Eigenvalues & \textbf{Ours} & Power Iter. & \textbf{Ours} & Power Iter. \\
		\midrule
		\textbf{1: Real Distinct Eigenvalues} & $-\lambda_k = \lambda_1.\; \lambda_k \in \mathbb{R}_{\leq0}$  & \textcolor{red}{$\checkmark$} & $\times$ & \textcolor{red}{$\checkmark$} & $\times$ \\
		\textbf{2: Complex Distinct Eigenvalues} & $|\lambda_k|  = \lambda_1.\; \lambda_k \in \mathbb{C}$ & \textcolor{red}{$\checkmark$} & $\times$ & \textcolor{red}{$\checkmark$} & $\times$ \\
		3: Repeated Complete Eigenvalues & $\lambda_k = \lambda_1.\; \lambda_k \in \mathbb{R}_{\geq0}$  & \textcolor{red}{$\checkmark$} & $\checkmark$ & $\times$ & $\times$ \\
		4: Repeated Defective Eigenvalues & $\lambda_k = \lambda_1.\; \lambda_k \in \mathbb{R}_{\geq0}$  & N.A. & N.A. & N.A. & N.A. \\
		\bottomrule
	\end{tabular}
	\label{tab:table_compare_PowerMethod}
\end{threeparttable}
\end{table}
In Case One where there are $N$ real and distinct eigenvalues,
power iteration is not guaranteed to converge
in particular when there is a $\lambda_k$ such that $\lambda_1 = - \lambda_k$.
Our method on the other hand is not affected by the presence of an opposite sign dominant eigenvalue,
which is the case for all bipartite graphs such as tree graphs.
In Case Two with complex distinct eigenvalues,
our method can still converge to the dominant eigenvalue and eigenvector pair when power iteration will fail.

Note that for Case Four in Table~\ref{tab:table_compare_PowerMethod},
we fill the row with N.A. or not applicable.
For power iteration,
the method requires the matrix to be diagonalizable or not defective.
Interestingly, the same constraint applies to our method as well when the dominant eigenvalue is defective.
This can be observed from the derivation in Equation~\ref{eq:ItdtvsIt} where we factor $I(t+\Delta t)$ with $I(t)$.
When we try to plug in the fastest growing term from Table~\ref{tab:table_fastestGrowingTerm},
the above derivation only works for Case One to Three,
i.e. when the matrix is diagonalizable.

\section{Experimental Results}

In this section,
we compare our method with power iteration method
experimentally.
We demonstrate our method's advantages for tree graphs, bipartite graphs, directed graphs with periods, and Markov chains with spider-traps.
For simplicity,
we use infection rate of $\beta=100\%$, time step-size of $\Delta t = 1$, and initial graph infection condition of $\mathbf{x}(0) = \bm{1}$.
Thus the Equation~\ref{eq:last_dance} to estimate the dominant eigenvalue for each secant line segment will simply be $\lambda_1 = e^{m} - 1$.
Same configurations are used for both our method and power iteration.
The number of iterations for our method is defined as the number of time steps used.
We measure the convergence of eigenvectors by calculating the angle between them based on their dot product.

\subsection{Tree Graph}

\begin{figure}[!ht]
	\centering
	\includegraphics[width=0.55\textwidth]{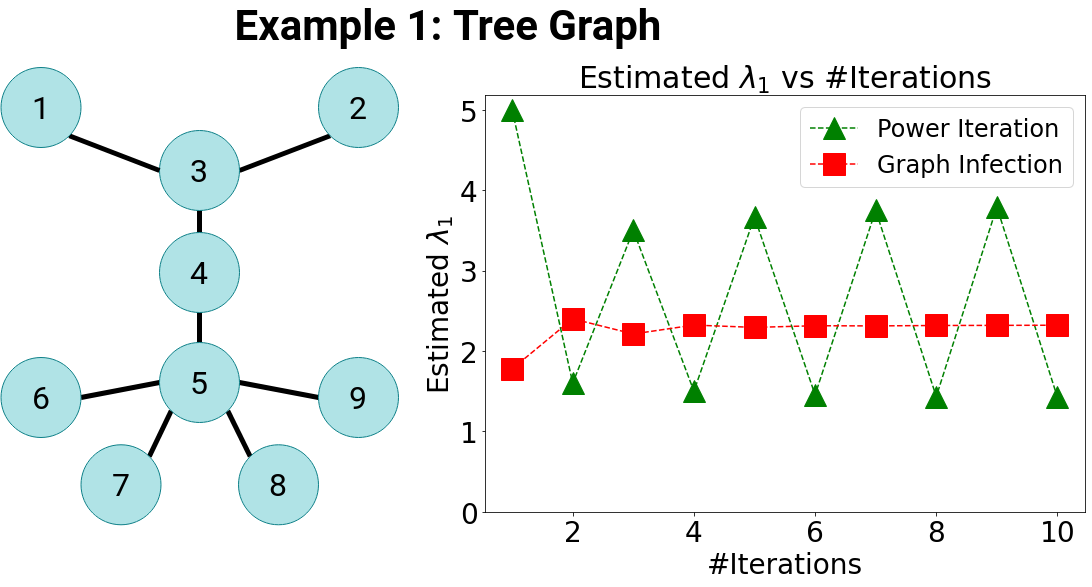}
	\caption{
	The dominant eigenvalue for this \textbf{tree graph} $\lambda_1 =\sqrt{4+\sqrt{2}}\approx2.327$.
	The estimated dominant eigenvalue by our method after ten iteration steps converges to $2.321$ approximately.
	Our estimated eigenvector also converges to ground truth eigenvector to near $0^{\circ}$ angle.
	Power iteration does not converge for either eigenvalue or eigenvector
	because of the presence of another eigenvalue of opposite sign $\lambda_2 = - \lambda_1$.}
	\label{fig:fig_results_tree}
\end{figure}

For unweighted tree graphs,
the eigenvalues of the adjacency matrix are all real since the matrix is symmetric.
Figure~\ref{fig:fig_results_tree} contains a tree graph with nine nodes.
The true dominant eigenvalue and eigenvector pair for this tree graph can be fairly estimated by our method after ten iteration steps.
However, since there always exists another eigenvalue of the same norm as $\lambda_1$ but of opposite sign,
power iteration does not converge to the true dominant eigenpair,
as shown in the figure's green plots.

\subsection{Bipartite Graph}

\begin{figure}[!t]
	\centering
	\includegraphics[width=0.55\textwidth]{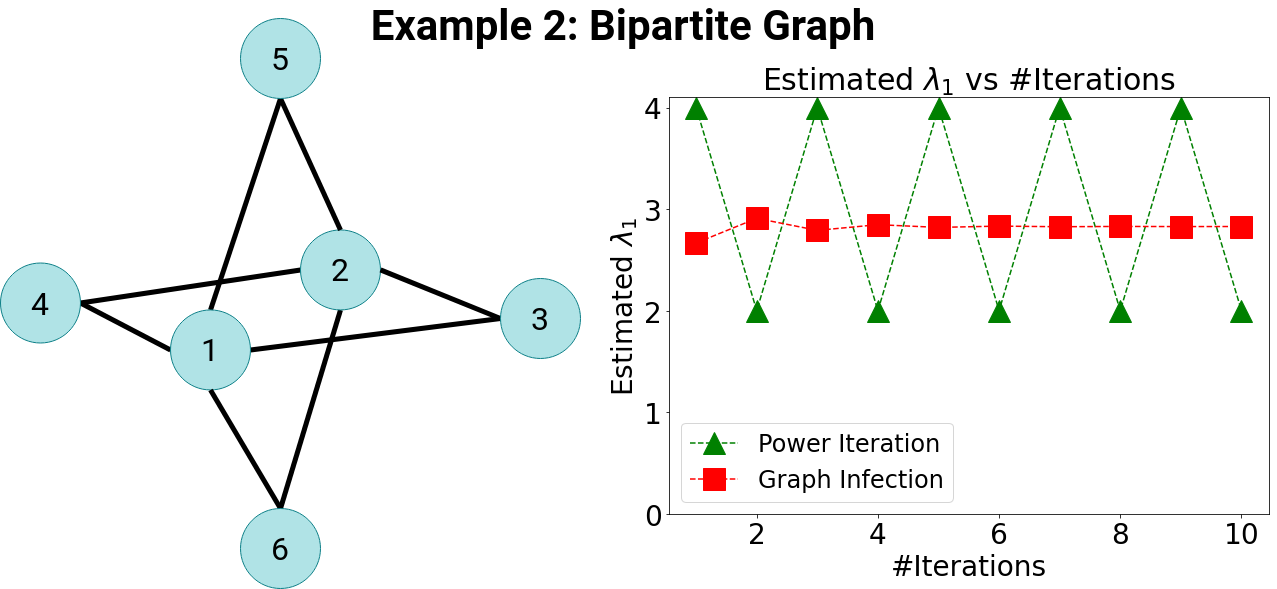}
	\caption{
	The true dominant eigenvalue for this \textbf{bipartite graph} $\lambda_1 =2\sqrt{2}\approx2.828$.
	The estimated dominant eigenvalue by our method after ten iteration steps converges to approximately $2.829$.
	Our estimated eigenvector also converges to ground truth within near $0^{\circ}$ angle.
	Power iteration does not converge for either eigenvalue or eigenvector
	because there is an eigenvalue of the same magnitude but of opposite sign $\lambda_2 = - \lambda_1$.}
	\label{fig:fig_results_bipartite}
\end{figure}

In fact, not just for tree graphs,
our method works better than power iteration for the superset of tree graphs: bipartite graphs.
Figure~\ref{fig:fig_results_bipartite} gives an example of a bipartite graph.
The bipartite graph with six nodes have a true dominant eigenvalue of around $2.83$,
which can be estimated by our method within ten iterations.
The associated eigenvector is also well estimated to within $0^{\circ}$ angle by dot product.
Compare that to the failed convergence pattern from power iteration method plotted in green.

\subsection{Directed Graph with Period of Two}

\begin{figure}
	\centering
	\includegraphics[width=0.56\textwidth]{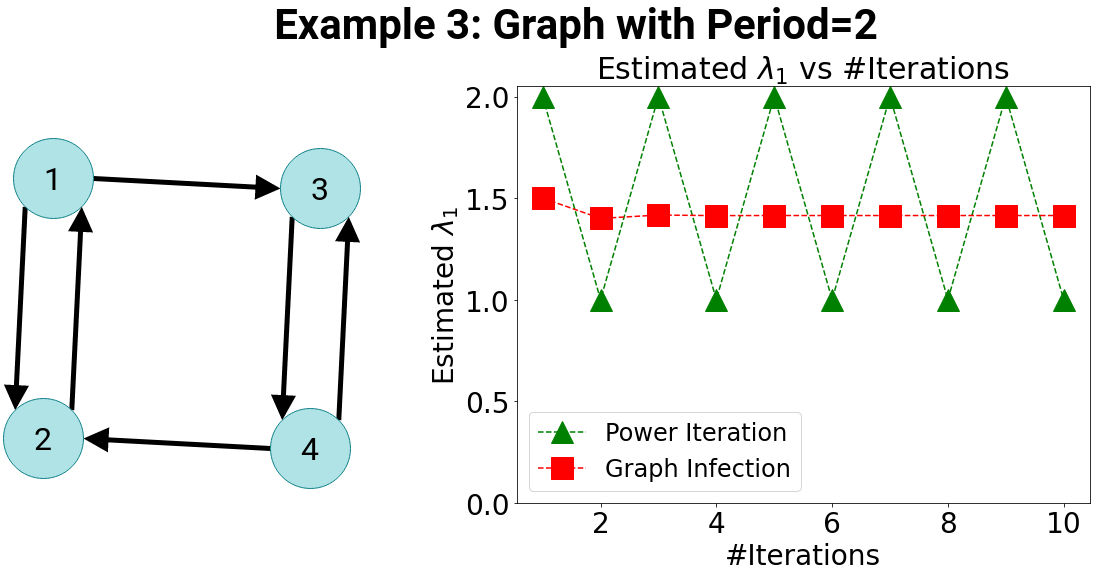}
	\caption{
	The true dominant eigenvalue for this \textbf{directed graph with period of two} $\lambda_1 =\sqrt{2}\approx1.414$.
	The estimated dominant eigenvalue by our method after ten iteration steps converges to around $1.414$.
	Our estimated eigenvector also converges to the true dominant eigenvector to within near $0^{\circ}$ angle.
	Power iteration does not converge for either eigenvalue or eigenvector
	because there are two eigenvalues of the same magnitude $\lambda_2 = - \lambda_1$.}
	\label{fig:fig_results_period2}
\end{figure}

Bipartite graphs can be abstracted as directed graphs with period of two when we convert all the edges into bi-directional edges.
Therefore we further experiment with the superset of bipartite graphs: directed graphs with period of two.
Figure~\ref{fig:fig_results_period2} shows one such graph.
This directed graph has a period of two,
and the true dominant eigenvalue $\lambda_1=\sqrt{2}$.
However, like bipartite graphs,
there exists another eigenvalue of the same magnitude but of opposite sign $\lambda_2 = - \sqrt{2}$.
Because of this, the power iteration cannot converge to the dominant eigenvalue nor the eigenvector.
Our method on the other hand can give a good estimate of the dominant eigenvalue and eigenvector pair within ten iterations.

\subsection{Markov Chain with Spider-trap}

Another prominent class of examples where our method outperforms the power iteration method is for transition probability graphs.
In particular, we look into Markov chains with `spider-traps',
a name coined to describe parts of the network from which a crawler cannot escape~\cite{rajaraman_ullman_2011}.
When a Markov chain contains a spider-trap,
the equilibrium distribution will be determined by the spider-trap in the long term.
Since there are as many eigenvalues equally spaced around the unit circle as the period of the Markov unichain~\cite{gallager2013stochastic_unichain_circle},
the period of the spider-trap will affect whether
the power iteration can converge.
If the period of the spider-trap is more than one,
then there are multiple eigenvalues of the same magnitude as the dominant eigenvalue $\lambda_1 = 1$.
Thus the power iteration method will not converge in Markov chains with periodic spider-traps.

Note that the dominant eigenvector $\bm{\mu}_1$ associated with the dominant eigenvalue $\lambda_1 = 1$ is an important attribute of the network that describes the equilibrium distribution of the Markov chain.

\subsubsection{Markov Chain with Period of Three}

\begin{figure}[!b]
	\centering
	\includegraphics[width=0.55\textwidth]{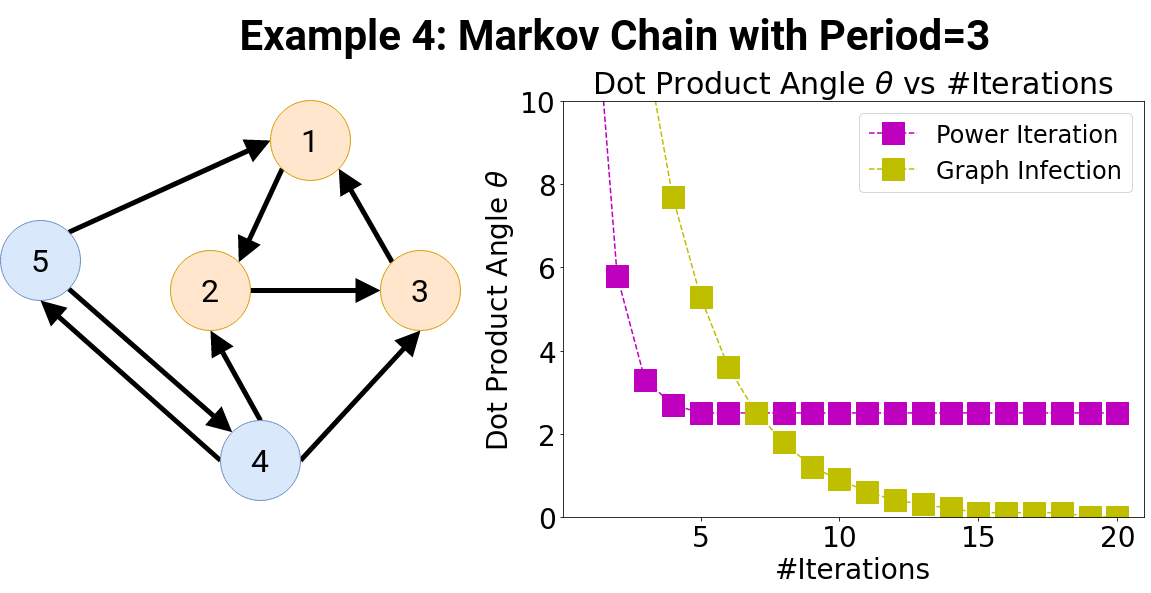}
	\caption{\textbf{Performance on Markov chain with period of three}.
	With our method,
	the dot product angle $\theta$ between the estimated dominant eigenvector and the true dominant eigenvector $\bm{\mu}_1$
	converges to almost $\theta\approx0^{\circ}$.
	Power iteration estimated dominant eigenvector fails to converge to true dominant eigenvector, with difference of $\theta\approx 2^{\circ}$,
	because there are two other complex conjugate eigenvalues of the same magnitude as $\lambda_1$.
    }
	\label{fig:fig_results_period3}
\end{figure}

For a Markov chain with period of three, see Figure~\ref{fig:fig_results_period3}.
There are five nodes in the Markov chain example.
Nodes 1-3 constitute a spider-trap with period of three.
This Markov chain has three eigenvalues of the same magnitude of 1, two of which are complex conjugates.
The remaining is the true dominant eigenvalue $\lambda_1 = 1$, and the associated dominant eigenvector $\bm{\mu}_1$ is the steady-state probability vector of the Markov chain.

Our previous discussion on Case Two of complex and distinct eigenvalues comes in handy here.
Because of the presence of complex eigenvalues of the same norm as the dominant real eigenvalue of $\lambda_1 = 1$,
power iteration method will not converge,
as shown in the magenta plots in Figure~\ref{fig:fig_results_period3}.
Our method does not have such issues and can estimate the dominant eigenvector $\bm{\mu}_1$ associated with $\lambda_1$ accurately as shown in yellow plots.

\subsubsection{Markov Chain with Period of Four}

\begin{figure}[!b]
	\centering
	\includegraphics[width=0.56\textwidth]{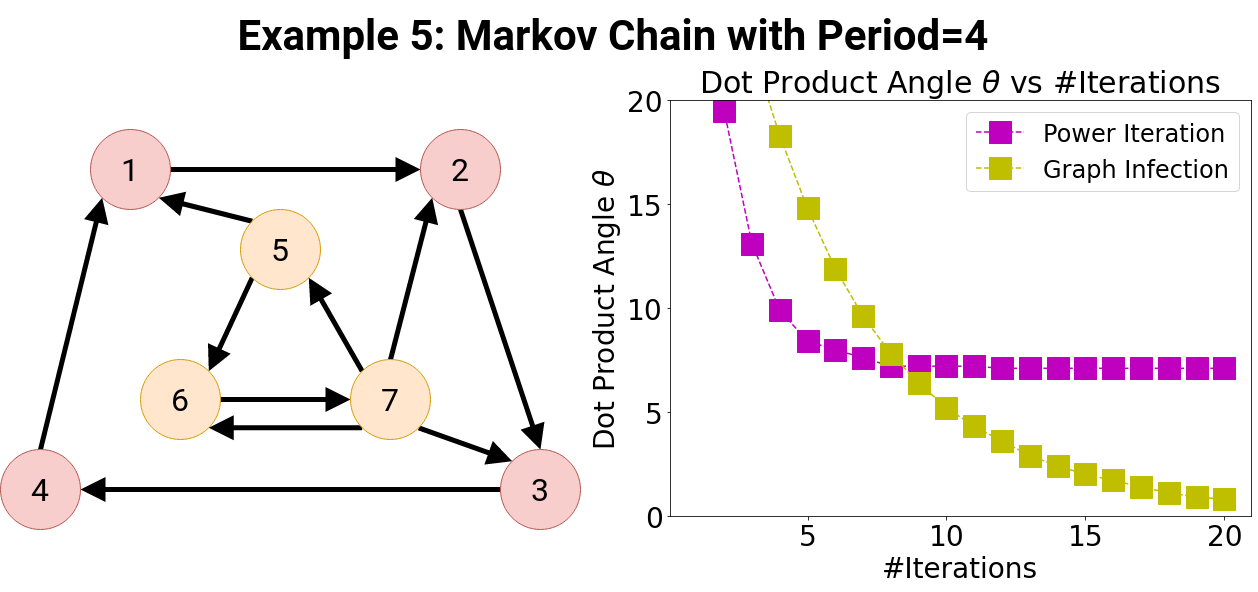}
	\caption{\textbf{Performance on Markov chain with period of four}.
	The dot product angle $\theta$ between the estimated dominant eigenvector and the true dominant eigenvector $\bm{\mu}_1$
	converges to almost $\theta\approx0^{\circ}$ with our method.
	Power iteration estimated dominant eigenvector fails to converge to true dominant eigenvector, with difference of $\theta\approx 7.1^{\circ}$,
	because there are four eigenvalues of the same magnitude of 1, $\lambda_1 = 1, \lambda_2 = -1, \lambda_3 = i, \lambda_4 = -i$.
    }
	\label{fig:fig_results_period4}
\end{figure}

There are seven nodes in the Markov chain depicted in
Figure~\ref{fig:fig_results_period4},
out of which Nodes 1-4 constitute a spider-trap with period of four.
For this Markov chain with a spider-trap of period four,
there are four eigenvalues of the same magnitude.
One of them is the $\lambda_1 = 1$ with the associated $\bm{\mu}_1$ steady-state vector.
The other three eigenvalues are spaced equally around the unit circle with values of $-1, i, -i$ respectively.

As discussed in the Case Two scenario,
the power iteration will not converge to the dominant eigenvector $\bm{\mu}_1$.
However, our method can estimate the eigenvector $\bm{\mu}_1$ as shown in the yellow vs magenta plots in Figure~\ref{fig:fig_results_period4}.

\section{Discussion}

There are a few key ideas that play important roles in our method.
One being the perspective of `reverse-engineering' the principal eigenvalue by solving a dynamical system based on graph infection epidemic model
(i.e. using infection dynamics to estimate a graph's eigenvalue)
.
Euler method is used as a viable way to iteratively approximate the solution to the matrix ODE that describes the dynamical system.
In particular, we apply the Perron-Frobenius theorem extended to non-negative matrices to derive the relationship between the non-negative real Perron-Frobenius eigenvalue and the fastest growing term of the ODE general solution.
Finally, it is interesting to see that when we plot the trends on total infection severity,
the slopes of secant line segments that arise from the Euler method
help us reveal the matrix's dominant eigenvalue and eigenvector pair.



\subsection{Relationship with NetworkX Eigenvector-centrality Implementation (Workaround on Power Iteration)}

NetworkX~\cite{NetworkX_SciPyProceedings_11} is a hugely influential Python package for network analysis.
In NetworkX, there is a built-in eigenvector centrality function that returns the dominant eigenvector of a graph.
This eigenvector-centrality feature has been available in NetworkX from as early as 2010 based on the NetworkX official GitHub repository git log.
However, until the release of NetworkX 2.0 beta 1 in August 2017,
the implementation of NetworkX's eigenvector-centrality algorithm has been the classical power iteration method which cannot converge when Case One or Case Two occurs, as we have discussed.

Interestingly, based on the GitHub discussions dating back to August 2015\footnote[1]{NetworkX GitHub issue \#1704: \url{https://github.com/networkx/networkx/issues/1704}} from their maintenance team of the NetworkX,
it was observed that addition of $\bm{A}$ by a positive multiple of the adjacency matrix,
such as
$(\bm{A} + \bm{I})$,
could alleviate the convergence problems for power iteration for graphs.
Since the NetworkX 2.0 beta 1 release in August 2017,
the eigenvector-centrality has thus been modified to be using power iteration on $(\bm{A} + \bm{I})$.

Here we want to highlight that NetworkX's workaround power iteration implementation can be thought of as a special case in our method.
Recall our Euler method and the finite difference approximation,
where we estimate the change of $\mathbf{x}(t)$ by
\begin{align*}
  \mathbf{x}(t + \Delta t) &\approx \mathbf{x}(t) + \frac{\mathrm{d} \mathbf{x}(t)}{\mathrm{d} t}\Delta t\\
  &= \mathbf{x}(t) + \beta \bm{A}\mathbf{x}(t)\Delta t \quad \text{by Equation~\ref{eq:Euler_dzdt_Az}}
\end{align*}
When we take time step size of 1 unit, $\Delta t = 1$,
infection rate to be $\beta = 100\%$,
the above finite difference equation becomes
\begin{align*}
  \mathbf{x}(t + 1) &\approx \mathbf{x}(t) + \frac{\mathrm{d} \mathbf{x}(t)}{\mathrm{d} t}\\
  &= \mathbf{x}(t) + \bm{A}\mathbf{x}(t)\\
  &= (\bm{A} + \bm{I})\mathbf{x}(t)
\end{align*}
In other words, NetworkX's workaround power iteration on $(\bm{A} + \bm{I})$ is in fact operating similar to a special context of our graph infection method.
Note that the consideration in NetworkX's shifted power iteration implementation is to shift the graph spectrum along the positive real axis direction.
Whereas our motivation stems from trying to solve an epidemic dynamical system iteratively.
Therefore, it is interesting that two conceptually different starting points eventually reached solutions that take a similar mathematical form!


\subsection{Matrix-Free Implementation}

We remark on the matrix-free nature of our method.
Implementation-wise, our method just sums all the node's infection severity after each time-step,
$I(t) = \sum_{i=1}^{N} x_i(t)$,
and computes the slope of the secant line $m$ by dividing the difference in log of total infection severity,
$\ln \left( I(t) \right)$,
over $\Delta t$.
Therefore our method is matrix-free, and does not require the whole network matrix to be stored inside the memory,
which is important for large graphs and can be parallelized.
Please refer to section Source Code~\ref{sec:source_code} for our example implementation, only requiring a few lines of code, in Python or R language.

\subsection{Extension to Dynamic Graphs}

It is possible to expand our scope to non-static graphs.
For example, there is a theoretical bound on the changes of the dominant eigenvector for strongly connected graphs under perturbation~\cite{dietzenbacher1988perturbations}.
Furthermore, Chen and Tong~\cite{chen2017eigen} proposed an algorithm that can effectively monitor how the dominant eigenvalue (computed by our method, say)
can be incrementally updated if the graph's $\bm{A}$ is perturbed,
which can be applicable for fast-changing large graphs.
More generally on the stability of eigenvalues for graphs,
Zhu and Wilson~\cite{zhu2005stability,wilson2008study} provide some
experiments on the stability of the spectrum to perturbations in the graph,
which can be useful heuristics in dynamic graphs.

\subsection{Limitations}

As can be seen from the row on Case Three in Table~\ref{tab:table_compare_PowerMethod},
when the dominant eigenvalue $\lambda_1$ is repeated and complete,
our method cannot guarantee the convergence to the associated dominant eigenvector $\bm{\mu}_1$.
But this is the same limitation faced by power iteration method for repeated complete dominant eigenvalue,
as the converged vector is spanned by the $p_k$ number of eigenvectors corresponding to the geometric multiplicity.

When the dominant eigenvalue is not only repeated but also defective,
as shown from the row on Case Four in Table~\ref{tab:table_compare_PowerMethod},
due to the complicated arrangements with generalized eigenvectors,
our derivation involving the fastest growing term and slope of secant line no longer holds when the dominant eigenvalue is defective.
Two things to note on the Case Four though.
First is that the same limitation imposed by defective dominant eigenvalue also applies to power iteration.
The Jordan form of the defective dominant eigenvalue also renders power iteration ineffective.
Secondly,
note that our method does not require the underlying matrix to be diagonalizable.
As long as the defective eigenvalue does not happen to be the dominant eigenvalue $\lambda_1$,
our method can converge to the dominant eigenvalue (Case Three),
and to its dominant eigenvector (Case One and Case Two).

\section{Conclusion}

We have proposed a novel method to estimate the dominant eigenvalue and eigenvector pair of any non-negative real matrix via graph infection.
To our knowledge, this is the first work that estimates dominant eigen-pair from the perspective of a dynamical system and matrix ODE.
Our method overcomes several limitations of the classical power iteration when the matrix has multiple dominant eigenvalues that are complex conjugates or are of opposite sign.
We believe our method can be adopted as an alternative to power iteration, especially for graphs.
It is our hope that this fresh perspective of `reverse-engineering' the dominant eigenvalue and eigenvector from matrix ODE of epidemic dynamical system can not only be of some practical use, but also to leave some food for thought in the eigenvalue algorithm literature.

\section*{Source Code}
\label{sec:source_code}
We provide example implementations of our method in Python or R language.
Source code is available at GitHub: \url{https://github.com/FeynmanDNA/Dominant_EigenPair_Est_Graph_Infection}.

\section*{Acknowledgement}
We thank Olafs Vandans, Chee Wei Tan, Johannes C. Paetzold, Houjing Huang, and Bjoern Menze for helpful comments and discussions.

%
%
%
\bibliographystyle{splncs04}
\bibliography{references}
%




\end{document}